\documentclass[a4paper]{article}
\usepackage[T1]{fontenc}
\usepackage[latin1]{inputenc}
\usepackage{amsmath, amsthm, amsfonts, amssymb, euscript, stmaryrd, graphicx, appendix}
\usepackage{color, pstricks, psfrag}
\usepackage{fullpage}
\usepackage{enumitem}

\def\independent{{\perp\!\!\!\!\perp}}
\def\simdist{\stackrel{\mathcal{L}}{\sim}}

\author{Sylvain \textsc{Corlay}\footnote{Bloomberg Quant Research, 731 Lexington Avenue, New York, NY 10022, USA. scorlay@bloomberg.net}}

\title{Properties of the Ornstein-Uhlenbeck bridge}
\date{January 21, 2014}
\newtheorem{theo}{Theorem}[section]
\newtheorem{prop}[theo]{Proposition}

\newcommand{\R}{\mathbb{R}}

\newcommand{\E}{\mathbb{E}}

\newcommand{\PP}{\mathbb{P}}

\newcommand{\var}{\operatorname{Var}}
\newcommand{\Proj}{\operatorname{Proj}}

\def\keywordname{{\bf Keywords:}} 
\newcommand{\keywords}[1]{\par\addvspace\baselineskip\noindent\keywordname\enspace\ignorespaces#1}

\graphicspath{{./}{figures/}}
\begin{document}
\maketitle
\begin{abstract}
\par This paper presents a study of the properties of the Ornstein-Uhlenbeck bridge, specifically, we derive its Karhunen-Loève expansion for any value of the initial variance and mean-reversion parameter (or mean-repulsion if negative). We also show that its canonical decomposition can be obtained using some techniques related to generalized bridges. Finally, we present an application to the optimal functional quantization of the Ornstein-Uhlenbeck bridge.
\end{abstract}
\section*{Introduction}
Let $\theta$ and $\mu$ be two real numbers and $\sigma > 0$. The Ornstein-Uhlenbeck process of long-term mean $\mu$, mean-reversion parameter $\theta$ and volatility $\sigma$ is defined as the solution of the S.D.E.
\begin{equation}\label{eq:solution}
d X_t = \theta (\mu - X_t) dt + \sigma dW_t,
\end{equation}
where $W$ is a standard Brownian motion and $X_0 \simdist \mathcal{N}\left(x_0, \sigma_0^2\right)$ is independent of $W$. We have
$$
X_t = X_0 e^{-\theta t} + \mu \left( 1 - e^{-\theta t}\right) + \int_0^t \sigma e^{\theta(u-t)} dW_u.
$$
\par \noindent For a finite horizon $T>0$, let us derive the covariance function of the Ornstein-Uhlenbeck bridge. For $\theta \neq 0$, we have
\begin{multline*}
\qquad \E\big[\big(X_t - \E[X_t | X_T] \big) \big(X_s - \E[X_s | X_T]\big) \big] = \frac{1}{2\theta} e^{-\theta (s+t)} \left(\left(2 \theta \sigma_0^2 - \sigma^2\right) + \sigma^2 e^{2\theta s \wedge t} \right) \\
- \frac{1}{2\theta} e^{-\theta T} \frac{\left( \sigma^2 e^{\theta t} + \left(2\theta \sigma_0^2 - \sigma^2 \right) e^{-\theta t}\right) \left( \sigma^2 e^{\theta s} + \left(2\theta \sigma_0^2 - \sigma^2 \right) e^{-\theta s}\right)}{\sigma^2 e^{\theta T} + \left(2\theta \sigma_0^2 - \sigma^2 \right) e^{-\theta T}},
\end{multline*}
and for $\theta = 0$, we have \quad $\E\big[\big(X_t - \E[X_t | X_T] \big) \big(X_s - \E[X_s | X_T]\big) \big] = \sigma_0^2 + \sigma^2 s \wedge t - \frac{\left(\sigma_0^2 + \sigma^2 t \right) \left(\sigma_0^2 +\sigma^2 s \right)}{\sigma_0^2 + \sigma^2 T}$.
\vspace{2mm}
\par While the Karhunen-Lo\`eve expansion of Brownian motion and the Brownian bridge are known, no closed-form expression is available in the general case of a Gaussian process with a continuous covariance function such as fractional Brownian motion. Asymptotics on the rates of decay of the eigenvalues can sometime be obtained as they are related with the mean-regularity of the process, the small-ball probabilities and the rate of decay of the quantization error. However, closed-form expressions are very useful for a variety of numerical applications such as the use of functional quantization for cubature. Examples of Gaussian processes for which a closed-form expression of the Karhunen-Lo\`eve expansion exists are available in \cite{Barczy} and \cite{DeheuvelsMartynov}. The case of the Ornstein-Uhlenbeck process is derived in \cite{CorlayPagesStratification}.
\par The paper is organized as follows: Section \ref{sec:semimartingale} covers background on the Ornstein-Uhlenbeck bridge and its canonical decomposition. Using a recent result on generalized bridges, we retrieve the canonical decomposition already derived by Barczy and Kern in \cite{BarczyKern2013}. In Section \ref{sec:karhunen-loeve}, we present the derivation of the Karhunen-Loève expansion of the Ornstein-Uhlenbeck bridge. In Section \ref{sec:functional_quantization} we derive the optimal functional quantization of Ornstein-Uhlenbeck bridges of various parameters. 
\section{Canonical decompositions in the enlarged filtration}\label{sec:semimartingale}
\subsection{Backgrounds on generalized bridges}
\par \noindent Let $(X_t)_{t \in [0,T]}$ be a continuous centered Gaussian semimartingale starting from $0$ on $(\Omega,\mathcal{A},\PP)$ and $\mathcal{F}^X$ its natural filtration. Fernique's theorem \cite{Fernique1970} ensures that $\int_0^T \E\left[X_t^2\right] dt < +\infty$.
\par \noindent As Alili in \cite{AliliGeneralizedBridge}, we are interested in the conditioning with respect to a finite family $\overline{Z}_T:=(Z_T^i)_{i \in I}$ of Gaussian random variables, which are the terminal values of processes of the form $Z_t^i = \int_0^t f_i(s) dX_s$, $i \in I$, for some finite set of bounded measurable functions $\overline{f} = (f_i)_{i \in I}$. A \emph{generalized bridge} for $(X_t)_{t \in [0,T]}$ corresponding to $\overline{f}$ with end-point $\overline{z} = (z_i)_{i\in I}$ is a process $\left(X_t^{\overline{f},\overline{z}}\right)_{t \in [0,T]}$ with distribution $X^{\overline{f},\overline{z}} \simdist \mathcal{L}\left(X \middle| Z_T^i = z_i, \ i\in I \right)$. 
\subsubsection{Gaussian semimartingales}\label{sec:semimartingale_generalized}
We denote by $\nu_{\left.\overline{Z}_T \middle| \left((X_t)_{t \in [0,s]}\right) \right.} : \mathcal{B}(\R^I) \times C^0\left([0,s],\R \right) \to \R_+,$ the transition kernel corresponding to the conditional distribution $\mathcal{L}\left(\overline{Z}_T \middle| \left((X_t)_{t \in [0,s]}\right)\right)$. We make the assumption ($\mathcal{H}$) that for every $s \in [0,T)$, this transition kernel is absolutely continuous with respect to the Lebesgue measure and we denote by $\Pi_{(x_u)_{u \in [0,s]},T}$ its density. This hypothesis is equivalent to assuming that the conditional covariance matrix $Q(s, T) = \E\left[\left(\overline{Z}_T - \E\left[\overline{Z}_T\middle|(X_u)_{u \in [0,s]}\right] \right)\left(\overline{Z}_T - \E\left[\overline{Z}_T \middle|(X_u)_{u \in [0,s]}\right] \right)^* \middle| (X_u)_{u \in [0,s]}\right]$ is invertible.
\begin{theo}[Radon-Nikodym derivative]\label{thm:alili}
\par \noindent Under the ($\mathcal{H}$)  hypothesis, for any $s \in [0,T)$, and for $\PP_{\overline{Z}_T}$-almost every $\overline{z}\in \R^I$, $\PP\left[\cdot \middle|\overline{Z}_T = \overline{z} \right]$ is equivalent to $\PP$ on $\mathcal{F}^X_s$ and its Radon-Nikodym density is given by
$$
\left.\frac{d\PP\left[\cdot \middle|\overline{Z}_T = \overline{z} \right]}{d\PP}\right|_{ \mathcal{F}^X_s} = \frac{\Pi_{(X_u)_{u\in [0,s]},T}(\overline{z})}{\Pi_{0,T}(\overline{z})}.
$$
\end{theo}
\begin{prop}[Generalized bridges as semimartingales]\label{prop:bridge_as_semimartingale}
\par \noindent Let us define the filtration $\mathcal{G}^{X,\overline{f}}$ by $\mathcal{G}^{X,\overline{f}}_t := \sigma\left(\overline{Z}_T, \mathcal{F}^X_t\right)$, the enlargement of the filtration $\mathcal{F}^X$ corresponding to the above conditioning. We consider the stochastic process $D_s^{\overline{z}} := \frac{d \PP\left[\cdot\middle|\overline{Z}_T = \overline{z} \right]}{d\PP}_{|\mathcal{F}_s^X} = \frac{\Pi_{(X_t)_{t\in [0,s]},T}(\overline{z})}{\Pi_{0,T}(\overline{z})}$ for $s \in [0,T)$. Then, under the ($\mathcal{H}$)  hypothesis, and the assumption that $D^{\overline{z}}$ is continuous, $X$ is a continuous $\mathcal{G}^{X,\overline{f}}$-semimartingale on $[0,T)$. 
\end{prop}
\par \noindent Theorem \ref{thm:alili} and Proposition \ref{prop:bridge_as_semimartingale} were first established in the Brownian case in \cite{AliliGeneralizedBridge} and extended to Gaussian semimartingales in \cite{CorlayPartialQuantization}. 
\subsubsection{Canonical decomposition}
\par Following the lines of \cite{CorlayPartialQuantization}, we define $L_t^{\overline{z}} := \int_0^t \frac{d\Pi_{(X_u)_{u\in [0,s]},T}(\overline{z})}{\Pi_{(X_u)_{u\in [0,s]},T}(\overline{z})}$. We have
$$
d \left\langle X, L^{\overline{z}} \right\rangle_s = d \left\langle X, \E \left[\overline{Z}_T \middle|(X_u)_{u \in[0,\cdot]} \right] \right\rangle_s Q(s,T)^{-1} \left(\overline{z}-\E \left[\overline{Z}_T \middle|(X_u)_{u \in[0,s]} \right] \right)^*. 
$$
Thus, if $X = V + M$ is the $(\mathcal{F}^X,\PP)$-canonical decomposition of $X$ then $M - \left\langle X, L^{\overline{z}} \right\rangle$ is a $\left(\mathcal{G}^{X,\overline{f}},\PP\left[\cdot \middle| \overline{Z}_T = \overline{z}\right] \right)$-martingale, and the canonical decomposition of $X$ in this filtration is given by $X = \left(V + \left\langle X, L^{\overline{z}} \right\rangle\right) + \left(M - \left\langle X, L^{\overline{z}} \right\rangle\right)$. In the case where $X$ is a Markov process, the expression for $\left\langle X, L^{\overline{z}} \right\rangle$ can be simplified: For every $j \in I$ there exists $g_j \in L^2 ([0,T])$ such that $\E \left[Z^j_T \middle|(X_u)_{u \in[0,s]} \right] = \int_0^s f_j(u) dX_u + g_j(s) X_s$. Hence, if one assumes that functions $(g_j)_{j \in I}$ have finite-variations, which is the case if $X$ is an Ornstein-Uhlenbeck process, then $d \left\langle X, \E \left[\overline{Z}_T \middle|(X_u)_{u \in[0,\cdot]} \right]\right\rangle_s = \left( \overline{f}(s) + \overline{g}(s) \right) d\langle X \rangle_s$, and thus
\begin{equation}\label{eq:markov_final_canonical_decomposition}
d \left\langle X, L^{\overline{z}} \right\rangle_s = \sum\limits_{i \in I} \left(f_i(s) + g_i(s) \right) \sum\limits_{j\in I} \left(Q(s,T)^{-1}\right)_{ij} \left(z_j - \E \left[Z^j_T \middle|(X_u)_{u \in[0,s]} \right] \right) d\langle X \rangle_s.
\end{equation}
\subsection{Centered Ornstein-Uhlenbeck bridge starting from $0$}
We perform the conditioning of a centered Ornstein-Uhlenbeck process $X$ starting from $0$ (meaning that $\mu = x_0 = \sigma_0 = 0$) by $X_T = z$. With the same notation as the previous section, we have $\E\left[X_T \middle| (X_u)_{u \in[0,s]} \right] = X_s e^{-\theta(T-s)}$ and thus $d \left\langle X, E\left[Z_T\middle|(X_u)_{u \in [0,\cdot]}\right] \right\rangle_s = e^{-\theta(T-s)}d\langle X \rangle_s = \sigma^2 e^{-\theta (T-s)}ds$. In this case, Equation \eqref{eq:markov_final_canonical_decomposition} simplifies to
\begin{equation}\label{eq:canon_interm1}
d\langle X,L^z \rangle_s = e^{-\theta (T-s)} \sigma^2 Q(s,T)^{-1} \left( \overline{z} - X_s e^{-\theta(T-s)} \right)ds.
\end{equation}
Moreover, $Q(s,T) = \E\bigg[\bigg(\int_s^T \sigma e^{\theta (u-T)}dW_u \bigg)^2\bigg] = \sigma^2e^{-2\theta T} \int_s^T e^{2 \theta u} du,$ and thus
$$
Q(s,T)^{-1} =
\left\{\begin{array}{cl}
\frac{2\theta}{\sigma^2} e^{-2 \theta T} \frac{1}{e^{2 \theta T} - e^{2 \theta s}} & \textnormal{if} \ \theta \neq 0,\\
\frac{1}{\sigma^2 (T - s)} & \textnormal{if} \ \theta = 0.\end{array}\right.
$$
\par \noindent If $\theta \neq 0$, plugging this into \eqref{eq:canon_interm1} yields
$d \langle X, L^z \rangle_s = 2 \theta \frac{z e^{\theta(T+s)} - X_s e^{2 \theta s}}{e^{2\theta T} - e^{2 \theta s}}ds$, and we finally obtain the canonical decomposition of Barczy and Kern \cite{BarczyKern2013},
\vspace{-5mm}
\begin{equation}\label{eq:canon_decomp_theta_non0}
dX_t = -\theta \coth(\theta (T-t)) X_t dt + \frac{\theta z}{\sinh(\theta (T-t))}dt + \overbrace{\left(\theta X_t dt  + dX_t + d\langle X,L^z \rangle_t \right)}^{\left( \mathcal{G}^W , \PP[\cdot|Z_T = z]\right) \textnormal{-martingale}}.
\end{equation}
By Lévy's characterization of Brownian motion, the martingale part is of the form $\sigma d \widetilde{W}_t$ where $\widetilde{W}$ is a $\left( \mathcal{G}^W , \PP[\cdot|Z_T =z]\right)$-Brownian motion. If $\theta = 0$, we retrieve the Brownian bridge $dX_t = \frac{z - X_t}{T - t} dt + \sigma d\widetilde{W}_t.$
\subsection{The case of a non-deterministic starting point $(\sigma_0^2 \neq 0)$}
\par The conditional distribution of $X_0$ knowing $X_T$, $\mathcal{N}\left(\E[X_0|X_T], \var(X_0 | X_T) \right)$ is completely determined by the covariance and expectation given in the introduction. We have $X_t = \left(X_0 e^{-\theta t} + \mu \left( 1 - e^{-\theta t}\right)\right) \quad \overset{\independent}{+} \widetilde{X}_t$, where $\widetilde{X}_t = \int_0^t \sigma e^{\theta(u-t)} dW_u$ is a centered Ornstein-Uhlenbeck process starting from $0$. $\widetilde{X}_t$ only depends on $X_T$ through its dependence on $\widetilde{X}_T = X_T - X_0 e^{-\theta T} - \mu \left( 1 - e^{-\theta t}\right)$. Thus, plugging $z = X_T - X_0 e^{-\theta T} - \mu \left( 1 - e^{-\theta t}\right)$ into \eqref{eq:canon_decomp_theta_non0} gives the  $\left( \mathcal{G}^W , \PP[\cdot|X_T]\right)$-canonical decomposition of $\widetilde{X}$. 
\section{Karhunen-Loève expansion}\label{sec:karhunen-loeve}
\par In this section, we derive the Karhunen-Loève expansion of the Ornstein-Uhlenbeck bridge of any initial variance or mean-reversion parameter. The method of derivation is the same as the one used for the Ornstein-Uhlenbeck process in \cite{CorlayPagesStratification}.
\par The covariance operator of the Ornstein-Uhlenbeck bridge is defined by $T^{OB} f \ (s) := \int_0^T c(s,t) f(t) dt$, where $c(s,t)$ is the covariance function $\E[(X_t - \E[X_t | X_T])(X_s - \E[X_s | X_T])]$ given in the introduction. 
\begin{prop}
If $f \in C([0,T])$ and $g := T^{OB} f$, then $g$ satisfies the boundary value problem
\begin{equation}\label{eq:boundary_value_problem}
\begin{array}{c}
g '' -\theta^2 g = -\sigma^2 f\\
g(T) = 0, \qquad \sigma_0^2 g'(0) = \left(\sigma^2 - \theta \sigma_0^2\right)g(0)\\
\end{array}
\end{equation}
Conversely, if $f \in C([0,T])$ and $g \in C^2([0,T])$ satisfy these three properties, then $g = T^{OB} f$.
\end{prop}
\par \noindent This is proved by differentiating twice under the integral signs and evaluating at $t=0$ and $t=T$.
\vspace{1mm}
\par \noindent As a consequence, the eigensystem $T^{OB}f = \lambda f$ amounts to solving $\lambda g'' + (\sigma^2 - \lambda \theta^2) g = 0$ with the same boundary conditions. Hence the Karhunen-Loève eigenvalues and unit eigenfunctions of the Ornstein-Uhlenbeck bridge are 
\begin{equation}\label{eq:OB_KL}
\lambda_n^{OB} = \frac{\sigma^2}{w_n^2 + \theta^2} \qquad \textnormal{and} \qquad e_n^{OB}(t) = \left(\frac{T}{2} - \frac{1}{4 w_n} \sin\left(2 w_n T\right)\right)^{-1/2} \sin(w_n (t-T)),
\end{equation}
where $(w_n)_{n \geq 1}$ are the strictly positive and increasingly sorted solutions to
\begin{equation}\label{eq:functional_eq}
(\sigma^2 - \theta \sigma_0^2) \sin(w T) = -w \sigma_0^2 \cos(w T).
\end{equation}
\begin{enumerate}
\item \textbf{Deterministic starting point $(\sigma_0 = 0)$}
\par Equation \eqref{eq:functional_eq} then amounts to $\sin(w T) = 0$, and thus $w_n = \frac{n \pi}{T}$ for $n \geq 1$. (This case has already been derived in \cite{DanilukMuchorskiKL}.)
\item \textbf{Non-deterministic starting point $(\sigma_0 \neq 0)$}
\begin{enumerate}
\item If $\sigma^2 = \theta \sigma_0^2$, Equation \eqref{eq:functional_eq} amounts to 
$
\cos(w T) = 0,
$
and thus the increasingly sorted positive solutions are $w_n = \frac{n \pi}{T} - \frac{\pi}{2 T}$ for $n \geq 1$. 
\item If $\theta \sigma_0^2 < \sigma^2$, Equation \eqref{eq:functional_eq} amounts to 
$
\underbrace{\left(\sigma^2 - \sigma_0^2 \theta \right)}_{>0}\tan(w T) = -\sigma_0^2 w,
$
\par \noindent and thus the increasingly sorted positive solutions satisfy $w_n \in \left] \frac{n \pi}{T} - \frac{\pi}{2T}, \frac{n \pi}{T} \right[$ for $n \geq 1$. 
\item If $\theta \sigma_0^2 > \sigma^2$, Equation \eqref{eq:functional_eq} amounts to 
$
\underbrace{\left(\sigma_0^2 \theta - \sigma^2\right)}_{>0}\tan(w T) = \sigma_0^2 w. 
$
\par \noindent There is a unique solution in each interval of the form $ \left] \frac{k \pi}{T}, \frac{k \pi}{T} + \frac{\pi}{2 T} \right[$ for $k \geq 1$. There is another solution on $\left]0, \frac{\pi}{2T} \right[$ if and only if $\sigma_0^2 > \sigma_0^2 \theta - \sigma^2$.
\end{enumerate}
\end{enumerate}
\par
\begin{figure}[!ht]
	\begin{center}
	\includegraphics[height=5cm]{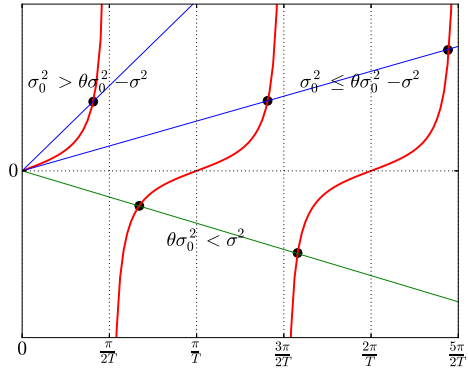}
 	\caption{Solutions to Equation \eqref{eq:functional_eq} in cases (b) and (c).}	\label{fig:tan_slope}
	\end{center}
\end{figure}

 \noindent In cases (b) and (c), the numerical value can then be computed using a root-finding method on the corresponding intervals. This is illustrated in Figure \ref{fig:tan_slope}. Using a certain rational approximation of $\tan$ proposed in \cite{CorlayPagesStratification}, a closed-form approximation of the solution is obtained as the root of a third-order polynomial and can be used as a starting point for the root-finding method. 
\vspace{-5mm}
\section{Functional quantization}\label{sec:functional_quantization}
\par The quantization of a random variable $X$ valued in a reflexive separable Banach space $(E,|\cdot|)$ consists in its approximation by $Y$ that takes finitely many values in $E$. We measure the resulting discretization error with the $L^2$ norm of the difference $|X-Y|$. If we settle on a fixed maximum cardinal $N$ for $Y(\Omega)$, the minimization of the error reduces to the optimization problem. 
\begin{equation}\label{eq:minimization_quantization}
\mathcal{E}_N(X,|\cdot|) = \min \left\{ \big\| \left|X-\Proj_\Gamma(X)\right| \big\|_2, \ \Gamma \subset E \ \textnormal{such that} \ |\Gamma| \leq N \right\},
\end{equation}
 A solution of (\ref{eq:minimization_quantization}) is an \emph{$L^2$-optimal quantizer} of $X$. 
\par Now let $X$ be a bi-measurable stochastic process on $[0,T]$ verifying $\int_0^T \E\left[|X_t|^2\right] dt < \infty$, which we see as a random variable valued in the Hilbert space $H = L^2([0,T])$. We assume that its covariance function $\Gamma^X$ is continuous. In the seminal paper \cite{LuschgyPagesFunctional3}, it is shown that, in the centered Gaussian case, linear subspaces $U$ of $H$ spanned by $N$-stationary quantizers correspond to principal components of $X$, in other words, are spanned by eigenvectors of the covariance operator of $X$, that is, its Karhunen-Loève eigenfunctions $\left(e_n^X\right)_{n \geq 1}$.
\par \noindent To perform optimal quantization, the Karhunen-Loève expansion is first truncated at a fixed order $m$ and then the $\R^m$-valued Gaussian vector constituted of the $m$ first coordinates of the process on its Karhunen-Loève decomposition is quantized. We have to determine the optimal rank of truncation $d^X(N)$ (the quantization dimension) and the optimal $d^X(N)$-dimensional quantizer of the first coordinates, $\bigotimes\limits_{j=1}^{d^X(N)} \mathcal{N}\left(0,\lambda^X_j\right)$. The minimal quadratic distortion $\mathcal{E}_N(X)$ is given by
$$
\mathcal{E}_N(X)^2 = \sum\limits_{j\geq m+1} \lambda_j^X + \mathcal{E}_N \bigg( \bigotimes\limits_{j=1}^m \mathcal{N}\left( 0, \lambda_j^X \right)\bigg)^2.
$$
\par
\begin{figure}[!ht]
	\begin{center}
	\includegraphics[height=5cm]{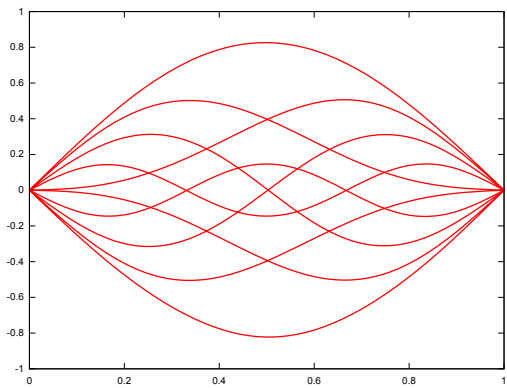}
 	\caption{Optimal quantization of the Ornstein-Uhlenbeck bridge with parameters $T = 1$, $\theta = 1$, $\sigma^2 = 1$, $\sigma^2_0 = 0$, $x_0 = \mu = X_T = 0$ and $N=10$.}	\label{fig:quantize}
	\end{center}
\end{figure}

\begin{figure}[!ht]
	\begin{center}
	\includegraphics[height=5cm]{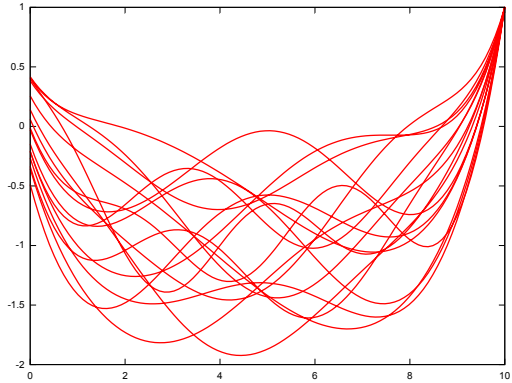}	
 	\caption{Optimal quantization of the Ornstein-Uhlenbeck bridge with parameters $T = 10$, $\theta = 1$, $\sigma^2 = 1$, $\sigma_0^2 = 1/2$, $x_0 = 0$,  $\mu = -1$, $X_T = 1$ and $N=16$.}\label{fig:quantize2}
	\end{center}
\end{figure}

 \noindent If the eigensystem $\left(e_n^X, \lambda_n^X \right)_{1 \leq n \leq d^X(N)}$ is known, we just need to perform the finite-dimensional quantization of $\bigotimes\limits_{j=1}^m \mathcal{N}\left(0, \lambda_j^X \right)$. Various algorithms have been devised to deal with this problem, among others,  Lloyd's algorithm \cite{LloydInformationTheory} and the Competitive Learning Vector Quantization (CLVQ)  \cite{PagesGaussianQuantization}. In Figures \ref{fig:quantize} and \ref{fig:quantize2}, we show optimal quantizers of the Orntein-Uhlenbeck bridge for different initial variances, mean-reversion parameters, volatilities and maturities.
\vspace{-4mm}
\par \noindent Moreover, using the rate of decay of the Karhunen-Lo\`eve eigenvalues of the Ornstein-Uhlenbeck bridge, and \cite[Theorem $2.2$]{LuschgyPagesFunctional2}, we see that the optimal quadratic quantization error of level $N$ of the Ornstein-Uhlenbeck bridge $\mathcal{E}_N(X)$ satisfies $\mathcal{E}_N(X) \sim K \log(N)^{-1/2}$ as $N \to \infty$, for some $K >0$.

\end{document}